\documentclass{article}

\usepackage{amsmath, amssymb, amsfonts, amsbsy, latexsym}
\usepackage{graphicx,stmaryrd}
\usepackage{epstopdf}

\newtheorem{theorem}{Theorem}[section]

\newtheorem{lemma}[theorem]{Lemma}
\newtheorem{proposition}[theorem]{Proposition}

\newtheorem{remark}[theorem]{Remark}
\newtheorem{example}[theorem]{Example}
\newtheorem{definition}[theorem]{Definition}

\def\C{\mathbb{C}}
\def\R{\mathbb{R}}

\def\<{\langle}
\def\>{\rangle}

\def\epsilon{\varepsilon}

\newcommand{\ip}[2]{\langle#1,#2\rangle}

\newcommand{\norm}[1]{\left\lVert#1\right\rVert}

\newcommand{\fc}{{\mathcal F}}

\newcommand{\hc}{{\mathcal H}}

\newcommand{\outp}[2]{{\llbracket #1,#2 \rrbracket }}
\newcommand{\ac}{\mathcal{A}}
\newcommand{\A}{\mathcal{A}}

\title{Phase Retrieval using Lipschitz Continuous Maps}

\author{Radu Balan$^{(1)}$ , Dongmian Zou$^{(2)}$\\
\hspace{-15mm}$^{(1)}$ Department of Mathematics and Center for Scientific Computation and Mathematical Modeling \\
 University of Maryland, College Park, MD 20742, USA \\
$^{(2)}$ Applied Mathematics, Applied Statistics and Scientific Computing Program \\
 University of Maryland, College Park, MD 20742, USA \\
}

\begin{document}

\maketitle

\begin{abstract}
In this note we prove that reconstruction from magnitudes of frame coefficients (the so called "phase retrieval problem") can be
performed using Lipschitz continuous maps. Specifically we show that when the nonlinear analysis map $\alpha:\hc\rightarrow\R^m$
 is injective, with $(\alpha(x))_k=|\ip{x}{f_k}|^2$, where $\{f_1,\ldots,f_m\}$ is a frame for the Hilbert space $\hc$, then
there exists a left inverse map $\omega:\R^m\rightarrow \hc$ that is Lipschitz continuous. Additionally we obtain the
Lipschitz constant of this inverse map in terms of the lower Lipschitz constant of $\alpha$. Surprisingly the increase in Lipschitz constant is independent of
the space dimension or frame redundancy.
\end{abstract}

\section{Introduction}\label{sec1}

Assume $\fc=\{f_1,f_2,\ldots,f_m\}$ is a frame (that is a spanning set) for the $n$-dimensional Hilbert space $H$. In this paper $H$ can be
a real or complex Hilbert space. The result applies to both cases, and the constants are the same.

Let $\alpha$ denote the nonlinear map
\begin{equation}\label{eq:1}
\alpha:H\rightarrow\R^m~~,~~\alpha(x)=\left(|\ip{x}{f_k}|^2\right)_{1\leq k\leq m}
\end{equation}
On $H$ we consider the equivalent replation $x\sim y$ iff there is a scalar $a$ of magnitude one, $|a|=1$, so that
$y=ax$. Let $\hat{H}=H/\sim$ denote the set of equivalence classes. Note $\hat{H}$ is equivalent to the cross-product between a real or complex projective
space ${\cal P}^{n-1}$ of dimension $n-1$ and the positive semiaxis $\R^{+}$.
The nonlinear map $\alpha$ is now well defined on $\hat{H}$
that, by abuse of notation, we denote also by $\alpha$. The {\em phase retrieval problem} (or the {\em phaseless reconstruction problem})
refers to analyzing when $\alpha$ is an injective map, and in this case to finding "good" left inverses.

The frame $\fc$ is said to be {\em phase retrievable} if the nonlinear map $\alpha$ is injective. In this paper we assume $\alpha$ is
injective (hence $\fc$ is phase retrievable). The problem is to extend a left inverse of $\alpha$ from $\alpha(\hat{H})$, the image of $\hat{H}$ through $\alpha$,
to the entire space $\R^m$ so that it remains Lipschitz continuous. A continuous map $f:(X,d_X)\rightarrow (Y,d_Y)$, defined between metric spaces $X$ and $Y$
with distances $d_X$ and $d_Y$ respectively, is Lipschitz continuous with Lipschitz constant $Lip(f)$ if
\[ Lip(f) := \sup_{x_1,x_2\in X} \frac{d_Y(f(x_1),f(x_2))}{d_X(x_1,x_2)} < \infty \]

Existing literature (e.g. \cite{BW13})  establishes that when the nonlinear map $\alpha$ is injective, it is also bi-Lipschitz for metric $d_1$ on $\hat{H}$ and Euclidian norm in $\R^m$.
Additionally, the nonlinear map $\sqrt{\alpha}:\hat{H}\rightarrow\R^m$, $\sqrt{\alpha}(x)=(|\ip{x}{f_k}|)_{1\leq k\leq m}$
is bi-Lipschitz in the case of real Hilbert space $H$ with distance $D_2$ on $\hat{H}$ and Euclidian distance on $\R^m$.
As a consequence of these results we obtain that a left inverse of $\alpha$ is Lipschitz when restricted to the image of $\hat{H}$ through $\alpha$.
In this paper we show that this left inverse of $\alpha$ admits a Lipschitz continuous extension to the entire $\R^m$. Surprisingly we obtain the Lipschitz constant
of this extension is just a small factor larger than the minimal Lipschitz constant, a factor that is independent of the 
dimension $n$ or the number of frame vectors $m$.

The organization of the paper is as follows. Section \ref{sec2} introduces notations and presents the main results. Section \ref{sec3} contains the proof
of these results.

\section{Notations. Statement of Main Results}\label{sec2}

The nonlinear map $\alpha$ naturally induces a linear map between the space $Sym(H)=\{T:H\rightarrow H~,~T=T^*\}$ of symmetric operators on $H$ and $\R^m$:
\begin{equation}
\label{eq:A}
\A:Sym(H)\rightarrow\R^m~~,~~\A(T)=(\ip{Tf_k}{f_k})_{1\leq k\leq m}
\end{equation}
This linear map has first been observed in \cite{BBCE07} and it has been exploited successfully in various paprs e.g. \cite{Bal2010,CSV12,Bal12a}.
Let $S^{p,}(H)$ denote the set of symmetric operators that have at most $p$ strictly positive eigenvalues and $q$ strictly negative eigenvalues.
In particular $S^{1,0}(H)$ denotes the set of non-negative symmetric operators of rank at most one:
\begin{equation} S^{1,0}(H) =\{ T\in Sym(H)~s.t.~\exists x\in H\, \forall y\in H~,~T(y)=\ip{y}{x}x \} \end{equation}
In \cite{Bal13a} we studied in more depth geometric and analytic properties of this set. In particular note
$\alpha(x)=\A(\outp{x}{x})$ where
\begin{equation} \outp{x}{y}=\frac{1}{2}(\ip{\cdot}{x}y+\ip{\cdot}{y}x)  \end{equation}
denotes the symmetric outer product between vectors $x$ and $y$.
The map $\alpha$ is injective if and only if $\A$ restricted to $S^{1,0}(H)$ is injective.

In previous papers \cite{Bal13a,BW13} we showed the following necessary and sufficient conditions for a frame to give phase retrieval.
\begin{theorem}\cite{Bal13a,BW13}
\label{equiv}
The following are equivalent:
\begin{enumerate}
\item The frame $\fc$ is phase retrievable;
\item $ker(\ac)\cap S^{1,1}(H) = \{0\}$;
\item There is a constant $a_0>0$ so that for every $u,v\in H$
\begin{equation}
\label{eq:lowerbound}
\sum_{k=1}^m \left|\frac{1}{2}( \ip{u}{f_k}|\ip{f_k}{v}+\ip{v}{f_k}\ip{f_k}{u}) \right|^2 \geq a_0 \left[ \norm{u}^2\norm{v}^2-(imag(\ip{u}{v}))^2 \right]
\end{equation}
\end{enumerate}
\end{theorem}

On the space $\hat{H}$ we consider two classes of metrics (distances) induced by corresponding distances on $H$ and $S^{1,0}(H)$ respectively.

The class of {\em vector space norm induced metrics}. For every $1\leq p \leq \infty$ and $x,y\in H$ define
\begin{equation}
D_p(\hat{x},\hat{y}) = \min_{|a|=1} \norm{x-ay}_p
\end{equation}
When no subscript is used, $\norm{\cdot}$ denotes the Euclidian norm, $\norm{\cdot}=\norm{\cdot}_2$.

The class of {\em matrix norm induced metrics}. For every $1\leq p\leq \infty$ and $x,y\in H$ define
\begin{equation}\label{eq:dp}
d_p(\hat{x},\hat{y}) = \norm{\outp{x}{x}-\outp{y}{y}}_p = \left\{ \begin{array}{rcl}
\mbox{$\left(\sum_{k=1}^n (\sigma_k)^p \right) ^{1/p}$} & for & \mbox{$1\leq p\leq \infty$} \\
\mbox{$\max_{1\leq k\leq n} \sigma_k $} & for & \mbox{$p=\infty$}
\end{array} \right.
\end{equation}
where $(\sigma_k)_{1\leq k\leq n}$ are the singular values of the matrix $\outp{x}{x}-\outp{y}{y}$, which is of rank at most 2.

Our choice in (\ref{eq:dp})
corresponds to the class of Schatten norms that extend to ideals of compact operators. In particular $p=\infty$ corresponds to the operator norm $\norm{\cdot}_{op}$ in
$Sym(H)$; $p=2$ corresponds to the Frobenius norm $\norm{\cdot}_{Fr}$ in $Sym(H)$; $p=1$ corresponds to the nuclear norm $\norm{\cdot}_1$ in $Sym(H)$:
\[ d_{\infty}(x,y) = \norm{\outp{x}{x}-\outp{y}{y}}_{op}~,~d_2(x,y)=\norm{\outp{x}{x}-\outp{y}{y}}_{Fr}~,~d_1(x,y)= \norm{\outp{x}{x}-\outp{y}{y}}_1 \]
Note the Frobenius norm $\norm{T}_{Fr} = \sqrt{trace(TT^*)}$ induces an Euclidian metric on $Sym(H)$. In \cite{Bal13a} Lemma 3.7 we computed eplicitely the eigenvalues of
 $\outp{x}{y}$. Based on these values, we can easily derive explicit expressions for these distances:
\[ d_{\infty}(\hat{x},\hat{y}) = \frac{1}{2}| \norm{x}^2-\norm{y}^2 | + \frac{1}{2}\sqrt{(\norm{x}^2+\norm{y}^2)^2-4|\ip{x}{y}|^2} \]
\[ d_{2}(x,y) = \sqrt{\norm{x}^4 + \norm{y}^4 - 2|\ip{x}{y}|^2 }\]
\[ d_1(x,y) = \sqrt{(\norm{x}^2+\norm{y}^2)^2-4|\ip{x}{y}|^2} \]
Since $\outp{x}{x}-\outp{y}{y}=\outp{u}{v}$ for $x=\frac{1}{2}(u+v)$ and $y=\frac{1}{2}(u-v)$, an equivalent condition to (\ref{eq:lowerbound}) is stated as follows:

\begin{enumerate}
\item[(iv)]\it There is a constant $a_0>0$ so that for every $x,y\in H$,
\begin{equation}
\label{eq:a0}
\norm{\alpha(x)-\alpha(y)} ^2 \geq a_0 (d_1 (x,y))^2
\end{equation}
\end{enumerate}

All metrics $D_p$ and $d_p$ induce the same topology as shown in the following result.
\begin{proposition}
\label{lem1}
1. For each $1\leq p\leq \infty$, $D_p$ and $d_p$ are metrics (distances) on $\hat{H}$.

2. $(D_p)_{1\leq p\leq\infty}$ are equivalent metrics, that is each $D_p$ induces the same topology on $\hat{H}$ as $D_1$.
Additionally, for every $1\leq p,q\leq\infty$
the embedding $i:(\hat{H},D_p)\rightarrow (\hat{H},D_q)$, $i(x)=x$,  is Lipschitz with Lipschitz constant
\begin{equation}\label{eq:L1}
L^D_{p,q,n}=max(1,n^{\frac{1}{q}-\frac{1}{p}}).
\end{equation}

3. For $1\leq p,q\leq\infty$, $(d_p)_{1\leq p\leq\infty}$ are equivalent metrics, that is each $d_p$ induces the same topology on $\hat{H}$ as $d_1$.
Additionally, for every $1\leq p,q\leq\infty$
the embedding $i:(\hat{H},d_p)\rightarrow (\hat{H},d_q)$, $i(x)=x$,  is Lipschitz with Lipschitz constant
\begin{equation}\label{eq:L2}
L^d_{p,q,n}=max(1,2^{\frac{1}{q}-\frac{1}{p}}).
\end{equation}

4. The metrics $D_p$ and $d_q$ are equivalent, that is they produce the same topology on $\hat{H}$. However neither the embedding
$i:(\hat{H},D_p)\rightarrow(\hat{H},d_q)$ nor $i:(\hat{H},d_q)\rightarrow(\hat{H},D_p)$ is Lipschitz, where $i(x)=x$.

5. The metric space $(\hat{H},d_p)$ is isometrically isomorphic to $S^{1,0}(H)$ endowed with Schatten norm $\norm{\cdot}_p$.
The isomorphism is given by the map
\begin{equation}\label{eq:kappa}
\kappa: \hat{H} \rightarrow S^{1,0}(H)~~,~~x\mapsto \outp{x}{x}.
\end{equation}
In particular the metric space $(\hat{H},d_1)$ is isometrically isomorphic to $S^{1,0}(H)$ endowed with the nuclear norm $\norm{\cdot}_1$.
\end{proposition}

\begin{remark}
\mbox{}

1. Note the Lipschitz bound $L^D_{p,q,n}$ is equal to the operator norm of the identity between $(\C^n,\norm{\cdot}_p)$ and $(\C^n,\norm{\cdot}_q)$:
 $L^D_{p,q,n}=\norm{I}_{l^p(\C^n)\rightarrow l^q(\C^n)}$.

2. Note the equality $L^d_{p,q,n}=L^D_{p,q,2}$.
\end{remark}

Theorem \ref{equiv} together with the previous proposition show that if the frame $\fc$ is phase retrievable then the nonlinear map (\ref{eq:1})
is bi-Lipschitz between metric spaces $(\hat{H},d_p)$ and $(\R^m,\norm{\cdot}_q)$. In particular, the Lipschitz constants between $(\hat{H},d_1)$ and
$(\R^m,\norm{\cdot}=\norm{\cdot}_2)$ are given by $\sqrt{a_0}$ and $\sqrt{b_0}$:
\begin{equation}
\sqrt{a_0} d_1(x,y) \leq \norm{\alpha(x)-\alpha(y)} \leq \sqrt{b_0} d_1(x,y)
\end{equation}
Clearly the inverse map defined on the range of $\alpha$  from metric space $(\alpha(\hat{H}),\norm{\cdot})$ to $(\hat{H},d_1)$:
\begin{equation}
\label{eq:inv1}
\tilde{\omega} : \alpha(\hat{H})\subset\R^m \rightarrow \hat{H}~~,~~\tilde{\omega}(c) = x~~{\rm if}~\alpha(x)=c
\end{equation}
 is Lipschitz with Lipschitz constant $\frac{1}{\sqrt{a_0}}$.
In this paper we prove that $\tilde{\omega}$ can be extended to the entire $\R^m$ as a Lipschitz map with Lipschitz constant that increases by a small factor.

The precise statement is given in the following Theorem which is the main result of this paper.

\begin{theorem}\label{maintheo}
Let $\fc=\{f_1,\ldots,f_m\}$ be a phase retrievable frame for the $n$ dimensional Hilbert space $H$, and let $\alpha:\hat{H}\rightarrow\R^m$ denote the
injective nonlinear analysis map $\alpha(x)=(|\ip{x}{f_k}|^2)_{1\leq k\leq m}$. Let $a_0$ denote the positive constant introduced in (\ref{eq:lowerbound}) and (\ref{eq:a0}).
Then there exists a Lipschitz continuous function $\omega:\R^m\rightarrow\hat{H}$ so that $\omega(\alpha(x))=x$ for all $x\in\hat{H}$.
For any $1\leq p,q\leq \infty$,  $\omega$ has an upper Lipschitz constant $Lip(\omega)_{p,q}$ between $(\R^m,\norm{\cdot}_p)$ and $(\hat{H},d_q)$ bounded by:
\begin{equation}\label{eq:Lpq}
Lip(\omega)_{p,q} \leq \left\{ \begin{array}{cc}
\mbox{$\frac{3+2\sqrt{2}}{\sqrt{a_0}} 2^{\frac{1}{q}-\frac{1}{2}}  max(1,m^{\frac{1}{2}-\frac{1}{p}})$} & \mbox{for $q\leq 2$}\\
\mbox{$\frac{3+2^{1+\frac{1}{q}}}{\sqrt{a_0}} max(1,m^{\frac{1}{2}-\frac{1}{p}})$} & \mbox{for $q> 2$}
\end{array} \right.
\end{equation}
Explicitly this means: for $q\leq 2$ and for all $c,d\in\R^m$:
\begin{equation}
\label{eq:Lippq}
d_q(\omega(c),\omega(d)) \leq \frac{3+2\sqrt{2}}{\sqrt{a_0}} 2^{\frac{1}{q}-\frac{1}{2}}  max(1,m^{\frac{1}{2}-\frac{1}{p}}) \norm{c-d}_p
\end{equation}
whereas for $q>2$ and for all $c,d\in\R^m$:
\begin{equation}
\label{eq:Lippq2}
d_q(\omega(c),\omega(d)) \leq \frac{1}{\sqrt{a_0}} (3+2^{1+\frac{1}{q}})  max(1,m^{\frac{1}{2}-\frac{1}{p}}) \norm{c-d}_p
\end{equation}
In particular, for $p=2$ and $q=1$ its Lipschitz constant $Lip(\omega)_{2,1}$
bounded by $\frac{4+3\sqrt{2}}{\sqrt{a_0}}=\frac{8.243}{\sqrt{a_0}}$:
\begin{equation}\label{eq:inv2}
d_1(\omega(c),\omega(d)) \leq \frac{4+3\sqrt{2}}{\sqrt{a_0}}\norm{c-d}
\end{equation}
\end{theorem}

The proof of Theorem \ref{maintheo}, presented in Section \ref{sec3}, requires construction of a special Lipschitz map. We believe this particular result is interesting in itself and
may be used in other constructions. Due to its importance we state it here:
\begin{lemma}
\label{prop:2}
Consider the spectral decomposition of any self-adjoint operator in $Sym(H)$, $A=\sum_{k=1}^d \lambda_{m(k)}P_k$, where
$\lambda_1\geq \lambda_2\geq\cdots\geq\lambda_n$ are the $n$ eigenvalues including multiplicities, and $P_1$,...,$P_d$ are the orthogonal projections
associated to the $d$ distinct eigenvalues. Additionally, $m(1)=1$ and $m(k+1)=m(k)+r(k)$, where $r(k)=rank(P_k)$ is the multiplicity of eigenvalue $\lambda_{m(k)}$.
Then the map
\begin{equation}
\label{eq:pi}
\pi: Sym(H) \rightarrow S^{1,0}(H)~~,~~\pi(A)=(\lambda_1-\lambda_2)P_1
\end{equation}
satisfies the following two properties: (1) for $1\leq p\leq \infty$, it is Lipschitz continuous from $(Sym(H),\norm{\cdot}_{p})$ to $(S^{1,0}(H),\norm{\cdot}_{p})$
with Lipschitz constant less than or equal to $3+2^{1+\frac{1}{p}}$;  (2) $\pi(A)=A$ for all $A\in S^{1,0}(H)$.
\end{lemma}

\begin{remark}
Numerical experiments suggest the Lipschitz constant of $\pi$ is smaller than 5 for $p=\infty$.
On the other hand it cannot be smaller than 2 as the following example shows.
\end{remark}

\begin{example}
If $A = \begin{pmatrix}
1 & 0 \\
0 & 1 \\\end{pmatrix}$,
$B = \begin{pmatrix}
2 & 0 \\
0 & 0 \\\end{pmatrix}$, then $\pi(A)=\begin{pmatrix}
0 & 0 \\
0 & 0 \\\end{pmatrix}$ and
$\pi(B) = \begin{pmatrix}
2 & 0 \\
0 & 0 \\\end{pmatrix}$.
Here we have $\norm{\pi(A)-\pi(B)}_{\infty}=2$ and $\norm{A-B}_{\infty}=1$. Thus for this example $\norm{\pi(A)-\pi(B)}_{\infty}=2\norm{A-B}_{\infty}$.
\end{example}

It is unlikely to obtain an isometric extension in Theorem \ref{maintheo}. Kirszbraun theorem \cite{WelWil75} gives a sufficient condition for
isometric extensions of Lipschitz maps.
The theorem states that isometric extensions are possible when the pair of metric spaces satisfy the Kirszbraun property, or the K property:
\begin{definition}
The Kirszbraun Property (K): Let $X$ and $Y$ be two metric spaces with metric $d_x$ and $d_y$ respectively. $(X,Y)$ is said to have Property (K) if for any pair of families of closed balls $\{B(x_i,r_i): i \in I\}$, $\{B(y_i,r_i): i \in I\}$, such that $d_y(y_i,y_j) \leqslant d_x(x_i,x_j)$ for each $i,j \in I$, it holds that $\bigcap B(x_i,r_i) \neq \emptyset \Rightarrow \bigcap B(y_i,r_i) \neq \emptyset$.
\end{definition}

If $(X,Y)$ has Property (K), then by Kirszbraun's Theorem we can extend a Lipschitz mapping defined on a subspace of $X$ to a Lipschitz mapping defined on $X$ while maintaining the Lipschitz constant. Unfortunately, if we consider $(X,d_X)=(\R^m,\norm{\cdot})$ and $Y=\hat{H}$, Property (K) does not hold for either $D_p$ or $d_p$.

\emph{Property (K) does not hold for $\hat{H}$ with norm $D_p$.} Specifically, $(\mathbb{R}^m,\mathbb{R}^n/\sim)$ does not have Property K.
\begin{example}
We give a counterexample for $m=n=2, p=2$: Let $\widetilde{y_1} = (3,1)$, $\widetilde{y_2} = (-1,1)$, $\widetilde{y_3} = (0,1)$ be the representatives of
three points $y_1$, $y_2$, $y_3$ in $\mathbb{R}^2/\sim$. Then $D_2(y_1, y_2) = 2\sqrt{2}$, $D_2(y_2, y_3) = 1$ and
$D_2(y_1, y_3) = 3$. Consider $x_1 = (0,0)$, $x_2 = (0,-2\sqrt{2})$, $x_3 = (-1,-2\sqrt{2})$ in $\mathbb{R}^2$ with the Euclidean distance,
then we have $\norm{x_1- x_2} = 2\sqrt{2}$, $\norm{x_2- x_3} = 1$ and $\norm{x_1- x_3} = 3$. For $r_1=\sqrt{6}$, $r_2=2-\sqrt{2}$, $r_3=\sqrt{6}-\sqrt{3}$,
we see that $(1-\sqrt{2},1+\sqrt{2}) \in \bigcap_{i=1}^{3}B(x_i,r_i)$ but $\bigcap_{i=1}^{3}B(y_i,r_i) = \emptyset$. To see $\bigcap_{i=1}^{3}B(y_i,r_i) = \emptyset$,
it suffices to look at the upper half plane in $\mathbb{R}^2$. If we look at the upper half plane $H$, then $B(y_1,r_1)$ becomes the union of two parts,
namely $B(\widetilde{y_1},r_1) \cup H$ and $B(-\widetilde{y_1},r_1) \cup H$, and $B(y_i,r_i)$ becomes $B(\widetilde{y_i},r_i)$ for $i=2$, $3$.
But $(B(\widetilde{y_1},r_1) \cup H) \cap B(\widetilde{y_2},r_2) = \emptyset$ and $(B(-\widetilde{y_1},r_1) \cup H) \cap B(\widetilde{y_3},r_3) = \emptyset$.
So we obtain that $\bigcap_{i=1}^{3}B(y_i,r_i) = \emptyset$.
\end{example}

\emph{Property (K) does not hold for $\hat{H}$ with norm $d_p$.} Specifically, $(\R^m,\mathbb{C}^n/\sim)$ does not have Property K.
\begin{example}
We give a counterexample for $m=4, n=2, p=2$: Let $\widetilde{y_1} = (1,1-i)$, $\widetilde{y_2} = (2,1-i)$ be the representatives of two points $y_1$, $y_2$ in $\mathbb{C}^2/\sim$. Then $d_2(y_1, y_2) = \sqrt{2}$. Consider $x_1 = (1,1,1,2)$ and $x_2 = (2,1,1,1)$ in $\mathbb{R}^4$ with the Euclidean distance, then $\norm{x_1- x_2} = \sqrt{2}$. For $r_1 = r_2 = 1/\sqrt{2}$, we see that $\bigcap_{i=1}^{2}B(x_i,r_i) = \{(\frac{3}{2},1,1,\frac{3}{2})\}$. On the other hand, $\bigcap_{i=1}^{2}B(y_i,r_i) = \emptyset$ since $\mathbb{C}^2/\sim$ can be isometrically embedded in $\mathbb{R}^4$ but $(\frac{3}{2},1,1,\frac{3}{2})$ is not in the map.
\end{example}

\begin{remark}
Using nonlinear functional analysis language (\cite{BenLin00}) Lemma \ref{prop:2} can be restated by saying that $S^{1,0}(H)$ is a 5-Lipschitz retract in $Sym(H)$.
\end{remark}

\section{Proofs of results}\label{sec3}
We start by proving proposition \ref{lem1}.

{\bf Proof of Proposition \ref{lem1}}

1. For $D_p$ obviously we have $D_p(\hat{x},\hat{y})\geqslant 0$ and $D_p(\hat{x},\hat{y})=0$ if and only if $\hat{x}=\hat{y}$. We also have $D_p(\hat{x},\hat{y})=D_p(\hat{y},\hat{x})$ since $\norm{x-ay}_{p}=\norm{y-a^{-1}x}_{p}$ for any $\hat{x}$, $\hat{y} \in H$, $|a|=1$. Also, for any $x$, $y$, $z \in H$, if $D_p(\hat{x},\hat{y})$ is achieved by $\norm{x-ay}_{p}$, $D_p(\hat{y},\hat{z})$ is achieved by $\norm{z-by}$, then $D_p(\hat{x},\hat{z}) \leqslant \norm{x-ab^{-1}z}_{p} = \norm{bx-az}_{p} \leqslant \norm{bx-aby}_{p} + \norm{aby-az}_{p} = D_p(\hat{x},\hat{y})+D_p(\hat{y},\hat{z})$. So $D_p$ is a metric.

Since $\norm{\cdot}_{p}$ in the definition of $d_p$ is the standard Schatten p-norm of a matrix, $d_p$ is also a metric.

2. For $p \leqslant q$, by H\"{o}lder's inequality we have for any $x=(x_1,x_2,...,x_n) \in H= \mathbb{C}^n$ that $\sum_{i=1}^{n}|x_i|^p \leqslant n^{(\frac{1}{p}-\frac{1}{q})}(\sum_{i=1}^{n}|x_i|^q)^{\frac{p}{q}}$. Thus $\norm{x}_{p} \leqslant n^{(\frac{1}{p}-\frac{1}{q})} \norm{x}_{q}$. Also since $\norm{\cdot}_{p}$ is homogeneous, if we assume $\norm{x}_{p}=1$ we have $\sum_{i=1}^{n}|x_i|^q \leqslant \sum_{i=1}^{n}|x_i|^p = 1$. Thus $\norm{x}_{q} \leqslant \norm{x}_{p}$. Therefore, we have $D_q(\hat{x},\hat{y})=\norm{x-a_1y}_{q} \geqslant n^{(\frac{1}{p}-\frac{1}{q})} \norm{x-a_1y}_{p} \geqslant n^{(\frac{1}{p}-\frac{1}{q})}D_p(\hat{x},\hat{y})$ and $D_p(\hat{x},\hat{y})=\norm{x-a_2y}_{p} \geqslant \norm{x-a_2y}_{q} \geqslant D_q(\hat{x},\hat{y})$ for some $a_1$, $a_2$ with magnitude $1$. Hence

\[D_q(\hat{x},\hat{y}) \leqslant D_p(\hat{x},\hat{y}) \leqslant n^{(\frac{1}{p}-\frac{1}{q})}D_q(\hat{x},\hat{y})\]

We see that $(D_p)_{1\leqslant p \leqslant \infty}$ are equivalent. The second part follows then immediately.

3. The proof is similar to 2. Note that there are at most 2 $\sigma_i$'s that are nonzero, so we have $2^{(\frac{1}{p}-\frac{1}{q})}$ instead of $n^{(\frac{1}{p}-\frac{1}{q})}$.

4. To prove that $D_p$ and $d_q$ are equivalent, we need only to show that each open ball with respect to $D_p$ contains an open ball with respect to $d_p$, and vise versa. By 2 and 3, it is sufficient to consider the case when $p=q=2$.

First, we fix $x \in H = \mathbb{C}^n$, $r > 0$. Let $R = \min(1,\frac{r}{(2\norm{x}_{\infty}+1)n^2})$. Then for any $\hat{y}$ such that $D_2(\hat{x},\hat{y})<R$, we take $y$ such that $\norm{x-y}<R$, then $\forall 1 \leqslant i,j \leqslant n$, $|x_i\overline{x_j}-y_i\overline{y_j}| = |x_i(\overline{x_j}-\overline{y_j})+(x_i-y_i)\overline{y_j}| < |x_i|R+R(|x_i|+R) = R(2|x_i|+R) \leqslant R(2|x_i|+1) \leqslant \frac{r}{n^2}$. Hence $d_2(\hat{x},\hat{y}) = \norm{xx^*-yy^*}_{2} < n^2 \cdot \frac{r}{n^2} = r$.

On the other hand, we fix $x \in H = \mathbb{C}^n$, $R>0$. Let $r = R^2/\sqrt{2}$. Then for any $\hat{y}$ such that $d_2(\hat{x},\hat{y}) < r$, we have
\begin{equation}
 (d_2(\hat{x},\hat{y}))^2 = \norm{x}^4+\norm{y}^4-2|\langle x,y \rangle|^2 < r^2 = \frac{R^4}{2}
\end{equation}

But we also have
\begin{equation}
 (D_2(\hat{x},\hat{y}))^2 = \min_{|a|=1} \norm{x-ay}^2 = \norm{x-\frac{\langle x,y \rangle}{|\langle x,y \rangle|}y}^2 = \norm{x}^2+\norm{y}^2-2|\langle x,y \rangle|
\end{equation}

So
\begin{equation}
 (D_2(\hat{x},\hat{y}))^4 = \norm{x}^4 + \norm{y}^4 + 2\norm{x}^2 \norm{y}^2 - 4(\norm{x}^2+\norm{y}^2)|\langle x,y \rangle|+4|\langle x,y \rangle|^2
\end{equation}

Since $|\langle x,y \rangle| \leqslant \norm{x}\norm{y} \leqslant (\norm{x}^2+\norm{y}^2)/2$, we can easily check that $(D_2(\hat{x},\hat{y}))^4 \leqslant 2(d_2(\hat{x},\hat{y}))^2 < R^4$. Hence $D_2(\hat{x},\hat{y}) < R$.

Thus $D_2$ and $d_2$ are indeed equivalent metrics. Therefore $D_p$ and $d_q$ are equivalent. Also, the imbedding $i$ is not Lipschitz: if we take $x=(x_1,0,\ldots,0) \in \mathbb{C}^n$, then $D_2(\hat{x},0) = |x_1|$, $d_2(\hat{x},0) = |x_1|^2$.

5. This follows directly from the construction of the map.

Q.E.D.
\vspace{5mm}

Next we prove Lemma \ref{prop:2}.

{\bf Proof of Lemma \ref{prop:2}}

(2) follows directly from the expression of $\pi$. We prove (1) below.

Let $A$, $B \in Sym(H)$ where $A=\sum_{k=1}^d \lambda_{m(k)}P_k$ and $B=\sum_{k'=1}^{d'} \mu_{m(k')}P_{k'}$. We now show that
\begin{equation}
\label{eq:lipc}
 \norm{\pi(A)-\pi(B)}_{p} \leqslant (3+2^{1+\frac{1}{p}})\norm{A-B}_{p}
\end{equation}

Assume $\lambda_1-\lambda_2\leq \mu_1-\mu_2$. Otherwise switch the notations for $A$ and $B$. If $\mu_1-\mu_2=0$ then $\pi(A)=\pi(B)=0$ and
the inequality (\ref{eq:lipc}) is satisfied. Assume now $\mu_1-\mu_2>0$. Thus $Q_1$ is of rank 1 and therefore $\norm{Q_1}_p=1$ for all $p$.
First note that
\begin{equation}
\label{eq:break}
 \pi(A)-\pi(B) = (\lambda_1-\lambda_2)P_1 - (\mu_1-\mu_2)Q_1 = (\lambda_1-\lambda_2)(P_1-Q_1) + (\lambda_1-\mu_1-(\lambda_2-\mu_2))Q_1
\end{equation}

Here $\norm{P_1}_{\infty} = \norm{Q_1}_{\infty} = 1$. Therefore we have $\norm{P_1-Q_1}_{\infty} \leqslant 1$ since $P_1$, $Q_1 \geqslant 0$. From that we have $\norm{P_1-Q_1}_p \leqslant 2^{\frac{1}{p}}$.

Also, by Weyl's inequality (see \cite{Bhatia} III.2)  we have $|\lambda_i-\mu_i| \leqslant \norm{A-B}_{\infty}$ for each $i$. Apply this to $i=1$, $2$ we get $|\lambda_1-\mu_1-(\lambda_2-\mu_2)| \leqslant |\lambda_1-\mu_1|+|\lambda_2-\mu_2| \leqslant 2\norm{A-B}_{\infty}$. Thus $|\lambda_1-\mu_1|+|\lambda_2-\mu_2| \leqslant 2\norm{A-B}_{\infty} \leqslant 2\norm{A-B}_p$.

Let $g:=\lambda_1-\lambda_2$, $\delta := \norm{A-B}_p$, then apply the above inequality to (\ref{eq:break}) we get
\begin{equation}
\label{eq:gp2d}
 \norm{\pi(A)-\pi(B)}_p \leqslant g\norm{P_1-Q_1}_p+2\delta \leqslant 2^{\frac{1}{p}}g+2\delta
\end{equation}

If $0 \leqslant g \leqslant (2+2^{-\frac{1}{p}})\delta$, then $\norm{\pi(A)-\pi(B)}_p \leqslant (2^{1+\frac{1}{p}}+3)\delta$ and we are done.

Now we consider the case where $g>(2+2^{-\frac{1}{p}})\delta$. We use holomorphic functional calculus and put
\begin{equation}
 P_1 = -\frac{1}{2\pi i}\oint_\gamma R_Adz
\end{equation}
and
\begin{equation}
 Q_1 = -\frac{1}{2\pi i}\oint_\gamma R_Bdz
\end{equation}
where $R_A = (A-zI)^{-1}$, $R_B = (B-zI)^{-1}$, and $\gamma = \gamma(t)$ is the contour given in the picture below and used also by \cite{ZB06}.

\begin{figure}[ht]
 \centering
 \includegraphics{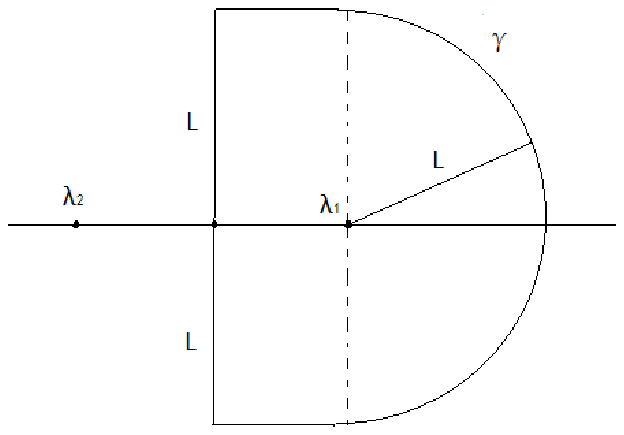}
\end{figure}

Therefore we have
\begin{equation}
\label{eq:pmq}
 \norm{P_1-Q_1}_p \leqslant \frac{1}{2\pi} \int_I \norm{(R_A-R_B)(\gamma(t))}_p |\gamma^{\prime} (t)| dt
\end{equation}

Now we have
\begin{equation}
\label{eq:ramrb}
 (R_A-R_B)(z) = R_A(z)-(I+R_A(z)(B-A))^{-1}R_A(z) = \sum_{n \geqslant 1} (-1)^n (R_A(z)(B-A))^n R_A(z)
\end{equation}
since for large $L$ we have $\norm{R_A(z)(B-A)}_{\infty} \leqslant \norm{R_A(z)}_{\infty}\norm{B-A}_p \leqslant \frac{\delta}{dist(z,\rho(A))} \leqslant \frac{2\delta}{g} < \frac{2}{2+2^{-\frac{1}{p}}} < 1$, where $\rho(A)$ denotes the spectrum of A.

Therefore we have
\begin{equation}
\begin{aligned}
\label{eq:nramrb}
 \norm{(R_A-R_B)(\gamma(t))}_p
 &\leqslant \sum_{n \geqslant 1} \norm{R_A(\gamma(t))}_{\infty}^{n+1} \norm{A-B}_p^n \\
 &=\frac{\norm{R_A(\gamma(t))}_{\infty}^2\norm{A-B}_p}{1-\norm{R_A(\gamma(t))}_{\infty}\norm{A-B}_p}
 <\frac{\norm{A-B}_p}{dist^2(\gamma(t),\rho(A))} \cdot (2^{1+\frac{1}{p}}+1)
\end{aligned}
\end{equation}
since $dist(\gamma(t),\rho(A)) \geqslant g/2$ for each t for large $L$. Here we used the fact that if we order the singular values of any matrix $X$ such that $\sigma_1(X) \geqslant \sigma_2(X) \geqslant \cdots$, then for any $i$ we have $\sigma_i(XY) \leqslant \sigma_1(X)\sigma_i(Y)$, and thus for two operators $X$, $Y \in Sym(H)$, we have $\norm{XY}_p \leqslant \norm{X}_{\infty} \norm{Y}_p$.

Hence by (\ref{eq:pmq}) and (\ref{eq:nramrb}) we have
\begin{equation}
\label{eq:npmq}
 \norm{P_1-Q_1}_p \leqslant (2^{\frac{1}{p}}+2^{-1})\frac{\norm{A-B}_p}{\pi} \int_I \frac{1}{dist^2(\gamma(t),\rho(A))}|\gamma^{\prime}(t)|dt
\end{equation}

By evaluating the integral and letting $L$ approach infinity for the contour, we have as in  \cite{ZB06}
\begin{equation}
\label{eq:atan}
 \int_I \frac{1}{dist^2(\gamma(t),\rho(A))}|\gamma^{\prime}(t)|dt = 2\int_0^{\infty}\frac{1}{t^2+(\frac{g}{2})^2}dt = \left[\frac{4}{g}\arctan\left(\frac{2t}{g}\right)\right]_0^{\infty} = \frac{2\pi}{g}
\end{equation}

Hence
\begin{equation}
\label{eq:npmq2}
 \norm{P_1-Q_1}_p \leqslant (2^{\frac{1}{p}}+2^{-1})\frac{\norm{A-B}_p}{\pi} \cdot \frac{2\pi}{g} = (2^{1+\frac{1}{p}}+1)\frac{\delta}{g}
\end{equation}

Thus by the first inequality in (\ref{eq:gp2d}) and (\ref{eq:npmq2}) we have $\norm{\pi(A)-\pi(B)}_{\infty} \leqslant
(3+2^{1+\frac{1}{p}})\delta$.

We have proved that $\norm{\pi(A)-\pi(B)}_{p} \leqslant (3+2^{1+\frac{1}{p}})\norm{A-B}_{p}$.
That is to say, $\pi:(Sym(H),\norm{\cdot}_p)\rightarrow(S^{1,0}(H),\norm{\cdot}_p)$
is Lipschitz continuous with Lipschitz constant less than or equal to $3+2^{1+\frac{1}{p}}$.

Q.E.D.
\vspace{5mm}

Now we are ready to prove Theorem \ref{maintheo}.

{\bf Proof of Theorem \ref{maintheo}}

We construct a map $\omega:(\R^m,\norm{\cdot}_p)\rightarrow (\hat{H},d_q)$ so that $\omega(\alpha(x))=x$ for all $x\in\hat{H}$, and $\omega$ is Lipschitz continuous.
We prove the Lipschitz bound (\ref{eq:Lpq}) which implies (\ref{eq:inv2}) for $p=2$ and $q=1$.

Set $M=\alpha(\hat{H})\subset\R^m$. By hypothesis, there is a map $\tilde{\omega}_1:M\rightarrow \hat{H}$ that is Lipschitz continuous and satisfies
$\tilde{\omega}_1(\alpha(x))=x$ for all $x\in\hat{H}$. Additionally, the Lipschitz bound between $(M,\norm{\cdot}_2)$ (that is, $M$ with Euclidian distance) and $(\hat{H},d_1)$
is given by $\frac{1}{\sqrt{a_0}}$.

First we change the metric on $\hat{H}$ from $d_1$ to $d_2$ and embed isometrically $\hat{H}$ into $Sym(H)$ with Frobenius norm (i.e. Euclidian metric):
\begin{equation}
\label{eq:p1}
 (M,\norm{\cdot}_2) \stackrel{\tilde{\omega}_1}{\longrightarrow} (\hat{H},d_1) \stackrel{i_{1,2}}{\longrightarrow} (\hat{H},d_2) \stackrel{\kappa}{\longrightarrow}
(Sym(H),\norm{\cdot}_{Fr})
\end{equation}
where $i_{1,2}(x)=x$ is the identity of $\hat{H}$ and $\kappa$ is the isometry (\ref{eq:kappa}).
We obtain a map $\tilde{\omega}_2:(M,\norm{\cdot}_2)\rightarrow (Sym(H),\norm{\cdot}_{Fr})$ of Lipschitz constant
$$Lip(\tilde{\omega}_2)\leq Lip(\tilde{\omega}_1)Lip(i_{1,2})Lip(\kappa^{-1}) = \frac{1}{\sqrt{a_0}},$$
 where we used $Lip(i_{1,2}) = L^d_{1,2,n}=1$ by (\ref{eq:L2}).

Kirszbraun Theorem \cite{WelWil75} extends isometrically $\tilde{\omega}_2$ from $M$ to the entire
$\R^m$ with Euclidian metric $\norm{\cdot}$. Thus we obtain a Lipschitz map $\omega_2:(\R^m,\norm{\cdot})\rightarrow (Sym(H),\norm{\cdot}_{Fr})$ of Lipschitz constant
$Lip(\omega_2)=Lip(\tilde{\omega}_2)\leq \frac{1}{\sqrt{a_0}}$ so that $\omega_2(\alpha(x))=\outp{x}{x}$ for all $x\in \hat{H}$.

The third step is to piece together $\omega_2$ with norm changing identities.

For $q\leq 2$ we consider the following maps:
\begin{equation}
\label{eq:p3}
(\R^m,\norm{\cdot}_p) \stackrel{j_{p,2}}{\longrightarrow}  (\R^m,\norm{\cdot}_2) \stackrel{{\omega}_2}{\longrightarrow} (Sym(H),\norm{\cdot}_{Fr})
\stackrel{\pi}{\longrightarrow} (S^{1,0}(H),\norm{\cdot}_{Fr}) \stackrel{\kappa^{-1}}{\longrightarrow}
(\hat{H},d_{2}) \stackrel{i_{2,q}}{\longrightarrow} (\hat{H},d_q)
 \end{equation}
where $j_{p,2}$ and $i_{2,q}$ are identity maps on the respective spaces that change the metric.
The map $\omega$ claimed by Theorem \ref{maintheo} is obtained by composing:
\[ \omega:(\R^m,\norm{\cdot}_p) \rightarrow (\hat{H},d_q)~~,~~\omega = i_{2,q}\cdot \kappa^{-1} \cdot \pi \cdot \omega_2 \cdot j_{p,2} \]
Its Lipschitz constant is bounded by
\[ Lip(\omega)_{p,q} \leq Lip(j_{p,2}) Lip(\omega_2)Lip(\pi)Lip(\kappa^{-1}) Lip(i_{2,q}) \leq max(1,m^{\frac{1}{2}-\frac{1}{p}}) \frac{1}{\sqrt{a_0}}\cdot (3+2\sqrt{2})\cdot 1
\cdot 2^{\frac{1}{q}-\frac{1}{2}} \]
Hence we obtained (\ref{eq:Lippq}). The other equation (\ref{eq:inv2}) follows for $p=2$ and $q=1$.

For $q>2$ we use:
\begin{equation}
\label{eq:p4}
(\R^m,\norm{\cdot}_p) \stackrel{j_{p,2}}{\longrightarrow}  (\R^m,\norm{\cdot}_2) \stackrel{{\omega}_2}{\longrightarrow} (Sym(H),\norm{\cdot}_{Fr})
\stackrel{I_{2,q}}{\longrightarrow} (Sym(H),\norm{\cdot}_{q})
\stackrel{\pi}{\longrightarrow} (S^{1,0}(H),\norm{\cdot}_{q}) \stackrel{\kappa^{-1}}{\longrightarrow}
(\hat{H},d_{q})
 \end{equation}
where $j_{p,2}$ and $I_{2,q}$ are identity maps on the respective spaces that change the metric.
The map $\omega$ claimed by Theorem \ref{maintheo} is obtained by composing:
\[ \omega:(\R^m,\norm{\cdot}_p) \rightarrow (\hat{H},d_q)~~,~~\omega = \kappa^{-1} \cdot \pi \cdot I_{2,q} \cdot \omega_2 \cdot j_{p,2} \]
Its Lipschitz constant is bounded by
\[ Lip(\omega)_{p,q} \leq Lip(j_{p,2}) Lip(\omega_2)Lip(I_{2,q})Lip(\pi)Lip(\kappa^{-1}) \leq max(1,m^{\frac{1}{2}-\frac{1}{p}}) \frac{1}{\sqrt{a_0}}\cdot
 1 \cdot (3+2^{1+\frac{1}{q}}) \cdot 1 \]
Hence we obtained (\ref{eq:Lippq2}) and this ends the proof. Q.E.D.

\section*{Acknowledgements}

The first author was supported in part by NSF grant DMS-1109498. He also acknowledges fruitful discussions with
Krzysztof Nowak and Hugo Woerdeman (both from Drexel University) who pointed out several references, with Stanislav Minsker (Duke University)
for pointing out \cite{ZB06} and \cite{DavKah69}, and Vern Paulsen (University of Houston) and Marcin Bownick (University of Oregon).


\begin{thebibliography}{1}  

\bibitem{ABFM12}
B.~Alexeev, A.~S.~Bandeira, M.~Fickus, D.~G.~Mixon, {\em Phase Retrieval with Polarization},
available online arXiv:1210.7752v1 [cs.IT] 29 Oct 2012.

\bibitem{App05}
D.~M. Appleby, \emph{Symmetric informationally complete-positive operator
  valued measures and the extended {C}lifford group}, J. Math. Phys.
  \textbf{46} (2005), no.~5, 052107, 29.

\bibitem{Bal99} R.~Balan, \emph{Equivalence relations and distances between Hilbert frames}, Proc. Amer.
Math. Soc. 127 (1999), no.~8, 2353–-2366.

\bibitem{Bal2009} R.~Balan, \emph{A Nonlinear Reconstruction Algorithm from Absolute Value
of Frame Coefficients for Low Redundancy Frames}, Proceedings of SampTA Conference, Marseille, France
May 2009.

\bibitem{Bal2010} R.~Balan, \emph{On Signal Reconstruction from Its Spectrogram}, Proceedings of the
CISS Conference, Princeton NJ, May 2010.

\bibitem{Bal12a} R.~Balan, \emph{Reconstruction of Signals from Magnitudes of Redundant Representations},
available online arXiv:1207.1134v1 [math.FA] 4 July 2012.

\bibitem{Bal13a} R.~Balan, \emph{Reconstruction of Signals from Magnitudes of Redundant Representations: The Complex Case},
available online arXiv:1304.1839v1 [math.FA] 6 April 2013.

\bibitem{BCE06}
R.~Balan, P.~Casazza, D.~Edidin, \emph{On signal reconstruction without phase},
 Appl.Comput.Harmon.Anal. \textbf{20} (2006), 345--356.

\bibitem{BCE07} R.~Balan, P.~Casazza, D.~Edidin,
\emph{Equivalence of Reconstruction from the Absolute Value of the
Frame Coefficients to a Sparse Representation Problem},
IEEE Signal.Proc.Letters, \textbf{14} (5)  (2007), 341--343.

\bibitem{BBCE07} R.~Balan, B.~Bodmann, P.~Casazza, D.~Edidin, \emph{Painless
reconstruction from Magnitudes of Frame Coefficients}, J.Fourier Anal.Applic., \textbf{15} (4) (2009), 488--501.

\bibitem{BW13} R.~Balan, Y.~Wang, \emph{Invertibility and Robustness of Phaseless Reconstruction},
available online arXiv:1308.4718v1.

\bibitem{BCMN13a}
A.S.~Bandeira, J.~Cahill, D.G.~Mixon, A.A.~Nelson, \emph{Saving phase: Injectivity and stability for phase retrieval},
available online arXiv:1302.4618v2 [math.FA] 14 Mar 2013.


\bibitem{BenLin00} Y.~Benyamini, J.~Lindenstrauss, {\bf Geometric Nonlinear Functional Analysis}, vol. 1, AMS Colloquium Publications, vol. 48, 2000.

\bibitem{Bhatia} R.~Bhatia, {\bf Matrix Analysis}, Graduate Texts in MAthematics 169, Springer-Verlag 1997.



\bibitem{BH13} B.~G.~Bodmann, N.~Hammen, \emph{Stable Phase Retrieval with Low-Redundancy Frames},
available online arXiv:1302.5487v1


\bibitem{CCPW13} J.~Cahill, P.G.~Casazza, J.~Peterson, L.~Woodland, \emph{Phase retrieval
by projections}, available online arXiv: 1305.6226v3


\bibitem{CSV12} E.~Cand\'{e}s, T.~Strohmer, V.~Voroninski, \emph{PhaseLift: Exact and Stable
Signal Recovery from Magnitude Measurements via Convex Programming},
Communications in Pure and Applied Mathematics vol. 66, 1241--1274 (2013).

\bibitem{CESV12} E.~Cand\'{e}s, Y.~Eldar, T.~Strohmer, V.~Voroninski, \emph{Phase
Retrieval via Matrix Completion Problem}, preprint 2011

\bibitem{Cass-artofframe}
P.~Casazza, \emph{The art of frame theory}, Taiwanese J. Math., (2) {\bf 4} (2000), 129--202.







\bibitem{Cah12} J.~Cahill, personal  communication, October 2012.

\bibitem{CamSei73}
P.~J. Cameron and J.~J. Seidel, \emph{Quadratic forms over {$GF(2)$}}, Indag. Math. \textbf{35} (1973), 1--8.

\bibitem{DavKah69} C.~Davis, W.M.~Kahan, \emph{Some new bounds on perturbation of subspaces},
 Bull. Amer. Math. Soc. vol. 75 (1969), no. 4, 863--868.



\bibitem{FMNW13} M.~Fickus, D.G.~Mixon, A.A.~Nelson, Y.~Wang, \emph{Phase retrieval from very few measurements},
available online arXiv:1307.7176v1

\bibitem{Fink04}
J.~Finkelstein, \emph{Pure-state informationally complete and ``really'' complete measurements},
Phys. Rev. A \textbf{70} (2004), no.~5, doi:10.1103/PhysRevA.70.052107


\bibitem{FP13} F.~Philipp, SPIE 2013 Conference Presentation, August 16, 2013, San Diego, CA.




\bibitem{HLO80}
M.~H.~Hayes, J.~S.~Lim, and A.~V.~Oppenheim, \emph{Signal Reconstruction
from Phase and Magnitude}, IEEE Trans. ASSP \textbf{28}, no.6 (1980),
 672--680.


\bibitem{HMW11}
T.~Heinosaari, L.~Mazzarella, M.~M.~Wolf, {\em Quantum Tomography under Prior Information},
arXiv:1109.5478v1 [quant-ph], 26 Sept 2011.

\bibitem{HG13} M.J.~Hirn, E.~Le~Gruyer, \emph{A general theorem of existence of quasi absolutely minimal Lipschitz extensions},
 arXiv:1211.5700v2 [math.FA], 8 Aug 2013.




\bibitem{Jm10} P.~Jaming, \emph{Uniqueness results for the phase retrieval problem of fractional Fourier transforms of variable
order}, preprint, arXiv:1009.3418.







\bibitem{MV13} D.~Mondragon, V.~Voroninski, \emph{Determination of all pure quantum states from a minimal number
of observables}, online arXiv:1306.1214v1 [math-ph] 5 June 2013.

\bibitem{Mil67}
R.~J.~Milgram, {\em Immersing Projective Spaces}, Annals of Mathematics, vol. \textbf{85},
no. 3 (1967), 473--482.

\bibitem{NQL82}
H.~Nawab, T.~F.~Quatieri, and J.~S.~Lim, \emph{Signal Reconstruction
from the Short-Time Fourier Transform Magnitude}, in Proceedings of
ICASSP 1984.





\bibitem{W12}
I.~Waldspurger, A.~d’Aspremont, S.~Mallat, \emph{Phase recovery, MaxCut and complex semidefinite programming},
Available online: arXiv:1206.0102



\bibitem{WelWil75} J.H.~Wells, L.R.~Williams, \textbf{Embeddings and Extensions in Analysis}, Ergebnisse der Mathematik und ihrer Grenzgebiete Band 84,
Springer-Verlag 1975.



\bibitem{ZB06} L.~Zwald, G.~Blanchard, \emph{On the convergence of eigenspaces in kernel Principal Component Analysis}, Proc. NIPS 05, vol. 18,
1649-1656, MIT Press, 2006.

\end{thebibliography}
\end{document}